\documentstyle[amscd]{amsart}
\newtheorem{thm}{Theorem}[section]
\newtheorem{lm}[thm]{Lemma}
\newtheorem{cor}[thm]{Corollary}
\newtheorem{pro}[thm]{Proposition}

\theoremstyle{definition}
\newtheorem{df}[thm]{Definition}

\theoremstyle{remark}
\numberwithin{equation}{section}

\def \R {\bold{R}}
\def \Z {\bold{Z}}

\def \J {\cal J}
\def \P {\cal P}
\def \CM {\cal M}
\def \hM {\hat{{\cal M}}}

\def \a {\alpha}

\def \g {\gamma}

\def \d {\delta}

\def \Lam {\Lambda}
\def \lam {\lambda}

\def \n {\nabla}
\def \o {\omega}
\def \p {\phi}
\def \r {\rho}
\def \s {\sigma}
\def \Si {\Sigma}
\def \t {\tau}
\def \T {\Theta}

\def \O {\Omega}

\def \bd {\partial}
\def \pj {\overline{\partial}_J}

\def \x {\times}
\def \ve {\varepsilon}

\begin{document}

\baselineskip.525cm

\title[Lagrangian embeddings]
{Lagrangian embeddings, Maslov indexes \\ and \\ 
integer graded symplectic Floer cohomology}

\author[Weiping Li]{Weiping Li}
\address{Department of Mathematics, Oklahoma State University \newline
\hspace*{.375in}Stillwater, Oklahoma 74078-0613}
\email{wli@@math.okstate.edu}
\dedicatory{Dedicated to R. Fintushel on the occasion of his fifty-th birthday}

\begin{abstract} We define an integer graded symplectic Floer cohomology
and a spectral sequence which are new invariants
for monotone Lagrangian sub-manifolds and exact isotopies. Such an integer
graded Floer cohomology is an integral lifting of the usual Floer-Oh
cohomology with $Z_{\Si (L)}$ grading. 
As one of applications of the spectral sequence, we offer an affirmative answer
to an Audin's question for oriented, embedded, monotone Lagrangian tori, i.e.
$\Si (L) = 2$.
\end{abstract}

\maketitle

\section{Introduction}

One of important problems in symplectic topology is to understand when a
manifold admits a Lagrangian embedding into $C^n$. M. Audin in \cite{au}
gave a series of obstructions to existence of such embeddings from purely
topological point of view. In \cite{au}, she asked whether
the following is true or not.
\medskip

\noindent{\bf Audin's Question}: Any embedded 
Lagrangian torus $L$ in $C^n$ has $\Si (L) =2$, where $\Si (L)$ is the 
positively minimal Maslov number of $L$, see \S 2 for definition.
\medskip

C. Viterbo in
\cite{vi2} derived a new obstruction to embedded Lagrangian tori, 
$2 \leq \Si (T^n) \leq n+ 1$. In particular, this answers
the Audin's question for $n =2$ since $\Si (L)$ is even for an
oriented Lagrangian
sub-manifold $L$. Viterbo used the symplectic action and 
Conley-Zehnder's finite
dimensional reduction technique to 
relate the Maslov index with the Morse index in \cite{vi}. His method for $T^n$
also works for compact manifold admitting a Riemannian metric with non-positive
sectional curvature.
In \cite{po1}, \cite{po2}, L. Polterovich obtained the same obstruction
on the Maslov index $1 \leq \Si (L ) \leq n+1$ for
more general manifold $L$ via Gromov's
pseudo-holomorphic curves.
\vspace{0.1in}

A new approach to the Lagrangian {\em rigidity} 
(relative to Lagrangian immersions,
Gromov's $h$-principle, {\em flexibility} \cite{alp}) is given by Y. G. Oh 
in \cite{oh1}. Using one of basic steps that Floer proved an 
Arnold conjecture for
monotone symplectic manifolds \cite{fl3}, Oh has constructed a {\em local}
symplectic Floer (co)homology such that 
there exists an approximation from local to
global symplectic Floer (co)homology by adding the holomorphic disks
or the quantum effect. This leads to 
a spectral sequence which converges to the usual symplectic Floer (co)homology. 
Thus the optimal rigidity for compact, embedded, monotone Lagrangian 
sub-manifolds is obtained 
in Theorem I \cite{oh1}, $1 \leq \Si (L) \leq n$. There are examples
in \cite{po2} which shows that the inequality $1 \leq \Si (L) \leq n$
can not be improved in general. Furthermore, Oh gives
a positive answer to the Audin's question for monotone Lagrangian tori of
dimension $\leq 24$ (Theorem III \cite{oh1}).
\medskip

From understanding Oh's result \cite{oh1} and \cite{oh2},
and Chekanov's result in \cite{ch}, we define an integer graded
symplectic Floer cohomology which can be thought as a combination
of Oh's and Chekanov's approaches.
In particular, we would like to replace the local symplectic Floer
cohomology by a more global symplectic invariant, and to relate the restricted
Floer (co)homology constructed in \cite{ch}.
This leads us to construct
an integer graded symplectic Floer cohomology in this paper
which is {\em global} and integral
lifting of the usual symplectic Floer cohomology, 
just like in the instanton theory 
developed by Fintushel and Stern \cite{fs2}. 
For special Lagrangian sub-manifolds
(representation varieties of handlebodies of an integral homology 3-sphere), 
Atiyah conjectured that the cohomologies from the instanton theory 
and the symplectic
theory are same. The Atiyah's conjecture is proved in \cite{ll} by R. Lee and 
the author. This leads us to borrow some ideas from the instanton Floer 
theory. Our construction of the integer graded Floer cohomology is related to
Chekanov's construction and Oh's construction via the Hofer's symplectic
energy (see \S 5). This may give an interesting way to understand 
the Hofer's symplectic energy through the integer graded Floer cohomology.
\medskip

Let $(P, \o)$ be a monotone symplectic manifold and $L$ is a monotone Lagrangian
sub-manifold in $(P, \o)$. Basically, we associate to $L$ a discrete set
$Im a_{\p}(Z_{\p})$, and for each $r \in {\R } \setminus Im a_{\p}(Z_{\p})
= {\R }_{L, \p}$ an $Z_2$-modules $I^{(r)}_*(L, \p; P)$ with a natural
integer grading. We also deal with more general case $(L_1, L_2)$ in \cite{ll2}
to get an integer graded symplectic Floer homology 
for a pair of monotone Lagrangian sub-manifolds. 
These $Z_2$-modules will depend on $r$ only through the
interval in ${\R }_{L, \p}$:

(i) $[r_0, r_1] \subset {\R }_{L, \p}$, then
$I^{(r_0)}_*(L, \p; P) = I^{(r_1)}_*(L, \p; P)$; 

(ii) $I^{(r)}_{* + \Si (L)}
(L, \p; P) = I^{(r+ \s (L))}_*(L, \p; P)$, where $\s (L) ( > 0)$ is the minimal 
number in $Im I_{\o }|_{\pi_2 (P, L)}$.
\medskip

\noindent{\bf Theorem A}. {\em For $\Si (L) \geq 3$, $L$ is a 
monotone Lagrangian sub-manifold in $(P, \o)$, 

(1) For any continuation $(J^{\lam }, {\p }^{\lam }) \in {\P }_1$ which
is regular at the ends, there exists an isomorphism
\[\p_{02}^n: I^{(r)}_n(L, \p^0; P, J^0) \to I^{(r)}_n(L, \p^1; P, J^1), \]
for all $n \in \Z$.

(2) There is a spectral sequence $(E^k_{n, j}, d^k)$ with
\[ E^1_{n,j}(L, \p; P, J) \cong I_j^{(r)}(L, \p; P, J), \ \ \
n\equiv j \ \ (mod \ \ \Si (L)), \]
and
\[E^{\infty }_{n,j}(L, \p; P, J) \cong F^{(r)}_n HF^j(L, \p; P, J)/
F^{(r)}_{n + \Si (L)}HF^j(L, \p; P, J) . \]

(3) The spectral sequence $(E^k_{n, j}, d^k)$ converges to the $Z_{\Si (L)}$
graded symplectic Floer cohomology $HF_*(L, \p; P, J)$, where
\[ d^k: E^k_{n, j}(L, \p; P, J) \to
E^k_{n + \Si (L) k +1, j +1}(L, \p; P, J) . \] 
}
\medskip

\noindent{\bf Theorem B}. {\em For each $k \geq 1$, $E^k_{n, j}(L, \p; P, J)$ are the symplectic invariant
under continuous deformations of $(J^{\lam }, {\p }^{\lam })$ within in the
set of continuations.}
\medskip

All the $E^k_{n, j}(L, \p; P, J) = E^k_{n, j}(L, \p; P)$,
$E^1_{n,j}(L,\p; P) = I_n^{(r)}(L,\p; P)$, for $k \geq 1, r \in {\R }_{L, \p}$
, are new symplectic invariants. They provided potentially interesting 
invariants for the symplectic topology of $L$. In particular the minimal $k$ for
which $E^k_{*, *} = E^{\infty }_{*, *}$ should be meaningful, denoted by
$k(L)$. 
Using the integer graded Floer cohomology, the spectral sequence,
and the new invariant $k(L)$, 
we give an affirmative
answer to Audin's question for monotone 
Lagrangian embedding torus in
$C^m$. 
\medskip

\noindent{\bf Theorem C}. {\em For an embedded, oriented, monotone Lagrangian
torus $L$ in $C^m$, we have
\[ \Si (L) = 2 . \] }

The question remains open for {\em non-monotone} Lagrangian torus
of dimension $\geq 3$. We use the full information about the spectral sequence
and some counting arguments to get the Theorem C unless 
$k(L) = \frac{m+1}{\Si (L)}$ (Theorem~\ref{unless}). Since the 
integer $\Si (L)$ is always even, so we first answer the Audin's question
for all even dimensional monotone embedded torus. Hence the complete proof
of Theorem C follows by studying for the odd 
dimensional case (Proposition~\ref{odd}).
Note that Oh in \cite{oh1} answered the Audin's question 
for monotone Lagrangian $L$ with $m \leq 24$ (Theorem III in \cite{oh1}). 
Our proof of Theorem C indicates the
interaction between $k(L)$ and $\Si (L)$ for a monotone Lagrangian sub-manifold
$L$.  
We also discuss the relation between the integer
graded symplectic Floer cohomology and the Chekanov's construction \cite{ch}.
Thus we obtain Chekanov's result by using the integer graded 
symplectic Floer cohomology. From these applications, 
the integer graded
symplectic Floer cohomology plays a uniform role in the work of Oh \cite{oh1}
and the work of Chekanov \cite{ch}. This is why we use $\s (L)$ to restrict
an energy band and $\Si (L) \geq 3$ to preserve the invariance from the
definition of $I^{(r)}_*(L, \p; P)$.
Our study suggests a possible relation between
the integer graded 
symplectic Floer cohomology and Hofer's symplectic energy for monotone
Lagrangian sub-manifolds. In fact we conjecture
that the Hofer's symplectic energy of a monotone Lagrangian sub-manifold $L$
with $\Si (L) \geq 3$ is a positive multiple of $\s (L)$ 
(More precisely, $e_H(L) = k(L) \s (L)$), 
where $\s (L)$
is the minimal symplectic action on $\pi_2(P, L)$. We will discuss this problem
elsewhere. It would be also
interesting to link the integer graded symplectic Floer cohomology with
the (modified) Floer cohomology with Novikov ring coefficients \cite{hs}, 
\cite{on}.
\medskip

The paper is organized as follows. In \S 2, we define the integer graded 
symplectic Floer cohomology for transversal Lagrangian intersections. Its invariant property under the symplectic continuations  is given in \S 3. Theorem A (1)
is proved in \S 3. Theorem A (2), (3) and Theorem B are proved in \S 4. In \S
5, we give some applications included Theorem C
and relations among the Hofer's energy, the Chekanov's construction and 
the Oh's results.

\section{Integer graded Floer cohomology for transversal Lagrangian
intersections}

Let $(P, \o)$ be a compact (or tamed) symplectic 
manifold with a closed non-degenerate
2-form $\o $. The 2-form $\o$ defines a second cohomology class 
$[\o ] \in H^2(P, {\R })$. By
choosing an almost complex structure $J$ on $(P, \o)$ such that
$\o ( \cdot , J \cdot )$ defines a Riemannian metric, we  have an
integer valued second cohomology class $c_1(P) \in H^2(P, Z)$ the first Chern
class. These two cohomology classes define two homomorphisms
\[I_{\o }: \pi_2(P) \to {\R }; \ \ \ \ \ \
I_{c_1}: \pi_2(P) \to {\Z }. \]
If $u : (D^2, \partial D^2) \to (P, L)$ is a smooth map of pairs, there is
a unique trivialization up to homotopy of the pull-back bundle $u^* TP
\cong D^2 \x C^n$ as a symplectic vector bundle. This trivialization
defines a map from $S^1 = \bd D^2$ to ${\Lam }(C^n)$ the set of
Lagrangians
in $C^n$. Let $\mu \in H^1({\Lam }(C^n), \Z)$ be the well-known Maslov
class. Then we define a map
\[ I_{\mu , L}: \pi_2(P, L) \to \Z, \]
by $I_{\mu , L}(u ) = \mu (\bd D^2)$, this Maslov index is invariant
under any symplectic isotopy of $P$.

\begin{df} (i) $(P, \o)$ is a monotone symplectic manifold if
\[ I_{\o} = \a I_{c_1}, \ \ \ \mbox{for some $\a \geq 0$}\]
(ii) A Lagrangian sub-manifold $L$ on $P$ is monotone if
\[ I_{\o } = \lam I_{\mu , L}, \ \ \ \mbox{for some $\lam \geq 0$}. \]
\end{df}

\noindent {\bf Remark:}
The monotonicity is preserved under the exact deformations of $L$.
By the canonical homomorphism $f: \pi_2(P) \to \pi_2(P,L)$,
one has
\[ I_{\o}(x) = I_{\o}(f(x)), \ \ \ \ \
I_{\mu , L}(f(x)) = 2 I_{c_1}(x), \]
where $x \neq 0 \in \pi_2(P)$. Therefore if $L$ is a
monotone
Lagrangian sub-manifold, then $P$ must be a monotone symplectic manifold and
$2 \lam =  \a$. In fact the constant $\lam$ does not depend
on  the Lagrangian $L$, but on $(P, \o)$ if $I_{\o}|_{
\pi_2(P)} \neq 0.$ For $\a = 0, \lam = 0$ cases, they are monotone in the
Floer's sense \cite{fl2}, \cite{fl3}. The definition of monotone Lagrangian
sub-manifolds is given by Oh in \cite{oh}.
\vspace{0.1in}

\begin{df} (i) 
Define $\s (L)$ be the positive minimal number in $Im I_{\o}|_{\pi_2(P, L)} 
\subset \R$;
Define $\Si (L)$ be the positive generator for the subgroup
$[\mu|_{\pi_2(P, L)}] = Im I_{\mu, L}$ in $\Z$. 

(ii) A Lagrangian sub-manifold
$L$ is called rational if $Im I_{\o}|_{\pi_2(P, L)} = \s (L) \Z$
is a discrete subgroup of $\R$ and $\s (L) > 0$ \cite{po4}. 
For a monotone Lagrangian, we
have $\s (L) = \lam \Si (L)$ for some $\lam > 0$.
\end{df}
 
Let $H: P \x \R \to \R$ be a smooth real valued function and let $X_H$ be 
defined by $\o (X_H, \cdot ) = dH$. Then the ordinary differential 
equation
\begin{equation} \label{hami} 
\frac{dx}{dt} = X_H(x(t)) ,
\end{equation} 
is called an Hamiltonian equation associated with the time-dependent
Hamiltonian
function $H$ or the Hamiltonian vector field $X_H$. It defines a family 
of diffeomorphisms of $P$ such that $x(t) = \p_{H,t}(x)$ solves 
equation (\ref{hami})
 for every $x$. A set ${\cal D}_{\o} = 
\{ \p_{H,1} | H \in C^{\infty}(P \x \R, 
\R)\}$ of all diffeomorphisms arising in this way is called the set of 
exact diffeomorphisms, an element in the set is called an exact isotopy.

\begin{df} For a given exact isotopy $\p = \{ \p_t\}_{0 \leq t \leq 1}$
on $(P, \o)$, we define
\[\O_{\p} = \{ z : I \to P \ \ \ | \ \ z(0) \in L, z(1) \in \p_1(L), 
[\p_t^{-1}z(t) ] = 0 \in \pi_1(P, L) \} . \]
\end{df}

\begin{pro} [Oh \cite{oh}]  \label{iden}
Let $L$ be monotone and $\p = \{ \p_t\}_{0 \leq t \leq 1}$
be an exact isotopy on $(P, \o)$. Let $u$ and $v$ be two maps from 
$[0,1] \x [0,1]$ to $\O_{\p}$ such that
\[u(\t , 0) , v(\t , 0) \in L, \ \ \ u(\t, 1 ) , v(\t, 1) \in \p(L) \]
\[u(0,t) = v(0,t) \equiv x, \ \ \ u(1,t) = v(1,t) \equiv y, \ \ \ 
x, y \in L \cap \p(L).\]
Then we have
\[I_{\o}(u) = I_{\o}(v) \ \ \mbox{if and only if} \ \ \ 
\mu_u(x,y) = \mu_v(x,y), \]
where $\mu_u$ is the Maslov-Viterbo index (see \cite{fl4} and \cite{vi}
for definition).

In particular, if $u, v$ are $J$-holomorphic curve with respect to an 
almost complex structure $J$ (may vary with time $t$) compatible with $\o$,
then 
\[ \int \|\n u \|^2_J = \int \| \n v\|^2_J \ \ \ \mbox{if and only if}
\ \ \ \mu_u(x,y) = \mu_v(x,y) . \]
\end{pro} 
This is Proposition 2.10 in \cite{oh}. Note that $\mu_u(x,y)$ is well-defined
mod $\Si (L)$ (see Lemma 4.7 \cite{oh}).
\vspace{0.1in}

The tangent space $T_z\O_{\p}$ consists of  vector fields $\xi $ of $P$ along 
$z$ which are tangent to $L$ at $0$ and to $\p(L)$ at $1$. Then $\o $ induces
a ``1-form'' on $\O_{\p}$.
\begin{equation} \label{1form}
 Da(z) \xi = \int_{S^1} \o (\frac{dz}{dt}, \xi(t)) dt .
\end{equation} 
This form is closed in the sense that it can be integrated locally to a real
function $a$ on $\O_{\p}$.
It  clearly vanishes
for all $\xi $ if and only if $z$ is a constant loop, i.e. $z(0)$ is 
a fixed point of $\p$; its critical points $Z_{\p }$ 
are the intersection points 
$L \cap \p(L)$. A critical point is non-degenerate if and only if 
the corresponding intersection is transversal.
\vspace{0.1in}

For a monotone Lagrangian sub-manifold $L$, an exact isotopy $\p$ and 
$k > 2/p$, consider the space of $L^p_k$-paths
\[{\P }^p_{k,loc}(L, \p;P) = \{ u \in L^p_{k, loc}(\T, P) \ \ 
| \ \ u({\R } \x \{0\}) \subset L, \ \  u({\R } \x \{1\} ) \subset \p(L) \}, \]
where ${\T }= {\R }\x [0,1] ={\R } \x iI \subset C$. 
Let $S_{\o}$ be the bundle of all 
$J \in End(TP)$ whose fiber is given by 
\[S_x = \{ J \in End(T_xP) \ \ | \ \ J^2 = - Id \ \ \mbox{and} \ \
\o(\cdot , J \cdot ) \ \ \mbox{is a Riemannian metric}\}, \]
and we denote the set of time-dependent almost complex structures by ${\J }= 
C^{\infty }([0,1] \x S_{\o})$. Define
\begin{equation}
\pj u(\t, t) = \frac{\bd u(\t, t)}{\bd \t} + J_t \frac{\bd u(\t,t)}{\bd t},
\end{equation} \label{holo}
on ${\P }^p_{k, loc}$ and then the equation $\pj u = 0$ is translationally
invariant in the variable $\t $. Denote ${\CM }= {\CM }_J(L, \p) = \{ u \in \O_{\p}
| \int_{\R \x I} | \frac{\bd u}{\bd \t}|^2 < \infty \}$ and
${\CM }(x, y) = \{ u \in {\CM }| \ \ \lim_{\t \to + \infty }u = x, \lim_{\t \to
- \infty }u = y; \ \ x, y \in Z_{\p} \}$.

${\CM }= \bigcup_{x, y \in L \cap \p(L)} {\CM }(x, y)$. If $L$ intersects $\p(L)$
transversely, then for each $x, y \in Z_{\p}$, there exist smooth Banach
manifold ${\P }(x, y) = {\P }^p_k(x,y) \subset {\P }^p_{k,loc}$ such that 
equation~\ref{1form} defines a smooth section $\pj$ of a smooth Banach space
bundle ${\cal L}$ over ${\P }(x,y)$ with fibers ${\cal L}_u = L^p_{k-1}(u^*TP)$,
and so that ${\CM }(x,y)$ is the zero set of $\pj $. The tangent space $T_u{\P}$
consists of all elements $\xi \in L^p_k(u^*(TP))$ so that $\xi (\t, 0) \in TL, 
\xi(\t, 1) \in T(\p (L))$ for all $\t \in \R$. The linearizations
\begin{equation}
E_u = D \pj (u): T_u {\P } \to {\cal L}_u
\end{equation}
are Fredholm operators for $u \in {\CM }(x,y)$. There is a dense set 
${\J }_{reg}(L,
\p(L)) \subset {\J }$ so that if $J \in {\J }_{reg}(L,\p(L)) $, then $E_u$ is 
surjective for all $u \in {\CM }(x, y)$. Moreover the Fredholm index of the 
linearization $Ind(E_u)$ is the same as the Maslov-Viterbo index $\mu_u(x,y)$.
In particular ${\CM }_J(x,y)$ becomes a smooth manifold with dimension equal to
$\mu_u(x,y)$ for $J \in {\J }_{reg}(L, \p(L))$.
See Proposition 2.1 \cite{fl2} and Theorem 1 \cite{fl4}.

\begin{thm}[Floer \cite{fl2} and Oh \cite{oh} Theorem 4.6] \label{focohomology}
Let $L$ be a monotone Lagrangian sub-manifold in $P$ and $\p = \{\p_t\}_{0 \leq
t \leq 1}$ be an exact isotopy such that $L$ is transversally intersects 
$\p_1(L)$. Let $C_*(L, \p; P, J)$ be the 
free generated $Z_2$-modules from $Z_{\p}$.
Suppose $\Si (L) \geq 3$. Then there exists a homomorphism
\begin{equation}
 \d: C_*(L, \p; P, J) \to C_*(L, \p; P, J)
\end{equation}
with $\d \circ \d = 0$ for $J \in {\J }_{reg}(L, \p(L))$. Define the 
Floer cohomology
$HF^*_J(L, \p; P, J)$ as the cohomology of $H^*(C_*(L, \p; P, J), \d )$,
a ${\Z }/ \Si(L)$-graded $Z_2$-modules. 
$HF^*_J(L, \p; P, J)$ is invariant under the continuation of $(J, \p)$,
denoted by $HF^{*}(L, \p; P)$ with $* \in Z_{\Si (L)}$.
\end{thm}

In order to extend this Floer-Oh theory to one with an integer grading we make 
use of an infinite cyclic cover $\tilde{\O}_{\p}$ of $\O_{\p}$ so that
the symplectic
action on $\tilde{\O}_{\p}$ and the Maslov index function on 
$\tilde{Z}_{\p}$ define as
$a: \tilde{\O}_{\p} \to \R$ and $\mu: \tilde{Z}_{\p} \to \Z$. From 
(\ref{1form}), a functional $a$ on ${\O }_{\p}$ is only defined up to $\s (L)Z$,
i.e. $a : {\O }_{\p} \to {\R }/ \s (L) Z$ from the different topology classes in
$\pi_2(P, L)$. The functional $a$ on ${\O }_{\p}$ and its lift on 
$\tilde{\O}_{\p}$ are clearly distinguished from its context. 

\begin{lm} \label{univ}
There exists a universal covering
 space $\tilde{\O}_{\p}$ of ${\O }_{\p}$ with transformation
group $\pi_2 (P, L)$.
\end{lm} 
Proof:  The function space ${\O }_{\p}$  has the homotopy type of a CW complex
and so an associated universal covering space. By Milnor's theorem 3.1 
in \cite{mi}, there is a universal covering space 
$\tilde{\O}_{\p}$ of ${\O }_{\p}$.
By definition of ${\O }_{\p}$, there is a homotopy $F(\t, t)$
of ${\p }^{-1}_t u(\t,t)$ to a constant
path in $(P, L)$ for any $u: I \to {\O }_{\p}$. Thus we can redefine the map
$u$ to yield a map $\overline{u}(\t, t) = u(2 \t, t)$ for $0 \leq \t \leq 1/2$;
and $\overline{u}(\t, t) = F(2 \t -1, t)$ for $1/2 \leq \t \leq 1$. Such a map
$\overline{u}: (D^2, \bd D^2) \to (P, L)$ defines an element in $\pi_2(P,L)$.
It is easy to check that $u \to \overline{u}$ is a bijective homomorphism
between $\pi_1(\O_{\p})$ and $\pi_2(P,L)$ (see also Proposition 2.3
in \cite{fl2}). \qed

Now the closed 1-form $Da(z)$ has a function up to a constant such that
$a : \tilde{\O}_{\p} \to \R$ is well-defined. Pick a point $z_0 \in L \cap 
\p(L)$ such that $a(z_0) = 0$ by adding a constant.
 For $g 
\in \pi_1(\O_{\p}) = \pi_2(P, L)$ we have 
\begin{equation} \label{deg}
a(g(x)) = a(x)
+ \deg (g) \s (L) , \end{equation}
where $\deg (g)$ is defined as $I_{\o}(g) = \deg (g) \s (L)$.
Let $Z_{\p} = \{ x \in L \cap \p(L) \ \  | \ \
[\p_t(x)] = 0 \in \pi_1(P,L) \}$. Let $Im a(Z_{\p})$ be the image of
$a$ of $Z_{\p}$; modulo $\s(L) \Z$, $Im a(Z_{\p})$ is a finite set. Thus
a set ${\R }_{L, \p} = {\R } \setminus Im a(Z_{\p})$ consists of the regular
values of the symplectic action $a$ on $\tilde{\O }_{\p}$. In this section, we
are going to construct an integer graded symplectic Floer cohomology for
every $r \in {\R }_{L, \p}$.
Given $x \in Z_{\p} \subset \O_{\p}$, let $x^{(r)} \in \tilde{Z}_{\p} \subset
\tilde{\O}_{\p}$ be the unique lift of $x$ such that $a(x^{(r)}) \in
(r, r + \s(L))$. Let $\mu^{(r)}(x) = \mu(x^{(r)}, z_0) \in \Z$ and define
the integral symplectic Floer cochain group
\begin{equation} \label{cohain}
 C^{(r)}_n(L, \p; P, J) = {\Z }_2 \{ x \in Z_{\p} \ \ \ | \ \ \ 
\mu^{(r)}(x) = n \in {\Z  }\} ,   
\end{equation}
as the free $Z_2$ module generated by $x \in Z_{\p}$ with the lift 
$x^{(r)}$ and $\mu(x^{(r)}, z_0) = n$.
If $\overline{z}_0$ is another choice of a based point and
$g(z_0) = \overline{z}_0$ for some covering transformation $g$, then the
corresponding choice of lift $\overline{x}^{(r)}$ of $x$ is just $g(x^{(r)})$.
Note that the integral Maslov index $\mu^{(r)}_u(x)$ is independent of the
choice of the based point $z_0$ used in the definition of $a$ by (\ref{deg}). 
The following 
lemma shows that the lift of the functional $a$ is compatible with a universal 
lift of ${\R }/ \s (L) Z$.
\begin{lm} \label{compatible}
The lift of the symplectic action over $\tilde{\O}_{\p}$ is 
compatible with the one of the Maslov index, i.e. for $g \in \pi_2(P, L)$
with $deg(g) = n$,
\[ a(g(z_0)) = n \s (L) \ \ \ \mbox{if and only if} \ \ \ \ 
\mu^{(r)}(g(z_0), z_0) = n \Si (L). \]
\end{lm} 
Proof: Let  $J$ be a compatible almost complex structure and
$\o (\cdot , J \cdot )$ be the corresponding Riemannian metric on $P$. Denote
$\n $ be the Levi-Civita connection of the metric $\o (\cdot , J \cdot )$.
Then $T_xL$ is an orthogonal complement of $JT_xL$. One can represent
$J_x$ to be a standard $J$ for suitable orthonormal basis in $T_xL$. Let $h$
be a parallel transport along the path $u(\t )$ ($u(\t, t)$ for each fixed
$t \in I$) in $\O_{\p}$. Then we obtain an isometry
\[h_{\t, t}: T_xP \to T_{u(\t, t)}P. \]
Define $J_{\t, t} = h_{\t, t}^{-1} \circ J_{u(\t, t)} \circ h_{\t, t}$. Then
we have a smooth map $f: I \x I \to SO(2m)$ such that $f_{\t, t}^{-1} \circ
J_{\t, t} \circ f_{\t, t} = J_x$. Set
\[ \tilde{L}(\t) = h_{\t, 0}^{-1}(T_{u(\t,0)}L), \ \ \
\tilde{\p(L)}(\t) = h_{\t, 1}^{-1}(T_{u(\t,1)}\p(L)). \]
Thus $f_{\t, 0}(\tilde{L}(\t)) = L(\t)$, $f_{\t, 1}(\tilde{\p(L)}(\t)) =
\p(L)(\t)$. The trivialization of $u^*TP$ by using the parallel transportation
$\{h_{\t, t}\}$ is given by \[u^*(T_xP) = I \x I \x T_xP = I \x I \x C^{m}.\]
Then there are two paths of Lagrangian subspaces $L(\t)$ and $\p(L)(\t)$ in
$T_xP = C^{m}$. Note that these two Lagrangian paths intersects transversally
at end points $\t = 0, 1$. There is a map $f_{\p}$
from $\O_{\p}$ to the space of pairs
of Lagrangian subspaces $\Lam (m) $
which has a universal covering $\tilde{\Lam }(m)$ (see \cite{ar}, 
\cite{clm} \S 1), where $f_{\p}(\{u(\t )\}) = \{L(\t), \p (L)(\t )\}, 0 \leq 
\t \leq 1$. 
For the map $f_{\p}$ from $\O_{\p}$ to $\Lam (m) $, there is a map from the
CW complex $\O_{\p}$ to $\Lam (m)$ from the obstruction theory. 
Hence there exists
a corresponding map $F$ between the universal covering space $\tilde{\O}_{\p}$
and the universal covering space $\tilde{\Lam }(m)$. From the choice of $z_0$,
we have $a(g(z_0)) = n \s (L)$. Here we have $u(0) = z_0, u(1) = g (z_0)$
and $\{u(\t )\}_{0 \leq \t \leq 1}$ corresponding to $g \in \pi_2 (P, L)$.

Note that we have $\deg : 
\pi_1 (\O_{\p}) \to \s (L) \Z$ and $Mas: \pi_1 (\Lam (m)) \cong \Si (L) \Z$
for our case. Hence there is an induced $g_1 \in \pi_1 (\Lam (m))$ such that 
$g_1 F = f_{\p } g$. $\mu^{(r)}(z_0, z_0) = 0$. We have the following
commutative diagram
\[ \begin{array}{ccc}
\pi_1 (\O_{\p }) & \stackrel{\pi_1(f_{\p })}{\longrightarrow} &
\pi_1 (\Lam (m)) \\
\downarrow  deg(g) & & \downarrow deg(g_1) \\
\s (L) {\Z } & \stackrel{F_*}{\longrightarrow} &
\Si (L) \Z .
\end{array} \]
So $I_{\o }(u_g) = n \s (L), I_{\mu , L}(u_g) = deg (g_1) \Si (L)$ by
the definitions of $\Si (L)$ and the Maslov index 
(see \cite{ar}). Thus the result follows from the
monotonicity and $\s (L) = \lam \Si (L)$.
\qed
 
\noindent{\bf Remark}: The index $\mu_u (x)$ depends on the trivialization
over $I \x I$, only the relative index does not depend on the trivialization.
So the choice of a single $z_0$ fixes the shifting in the integer graded Floer
cochain complex.
 
\begin{df} The integral Floer coboundary map 
$\bd^{(r)}: C^{(r)}_{n-1}(L, \p; P, J) \to C^{(r)}_{n}(L, \p; P, J)$ 
is defined by
\[ \bd^{(r)}x = \sum_{y \in C^{(r)}_{n}(L, \p; P, J)} \# \hM (x, y) \cdot
 y, \]
where ${\CM }(x,y)$ denote the union of the components of 1-dimensional moduli 
space of $J$-holomorphic curves and ${\hM }(x,y) = {\CM }(x, y)/{\R }$ is a zero
dimensional moduli space modulo $\t$-translationally invariant. 
$\# {\hM }(x, y)$
counts the points modulo 2.
\end{df} \label{integ}

Note that the coboundary map $\bd^{(r)}$ only counts part of the Floer-Oh's 
coboundary map in \cite{oh}. We are going to show that 
$\bd^{(r)} \circ \bd^{(r)} = 0$.
The corresponding cohomology groups are the integer graded symplectic
Floer cohomology $I^{(r)}_*(L, \p; P, J), * \in \Z$.

\begin{pro} \label{mas}
If $u \in {\P }(x, y)$ for $x, y \in Z_{\p}$ and $\tilde{u}$ is any 
lift of $u$, then $\mu_{\tilde{u}}= \mu^{(r)}(y) - \mu^{(r)}(x)$.
\end{pro}
This follows from Lemma~\ref{compatible} or Proposition 2.4 \cite{fl2} 
and Proposition 2.10 \cite{oh}. 
\begin{pro}[Oh , Proposition 4.1 and 4.3 \cite{oh}] \label{bub} 
Suppose that two Lagrangian sub-manifolds 
$L$ and  $\p(L)$ intersect transversally
and $\Si (L) \geq 3$. Then there is a dense subset 
${\J }_*(L, \p) \subset {\J }_{reg} (L, \p)$ of $\J$ such that (1) the zero 
dimensional component of ${\hM }(x, y)$ is compact and (2) the one dimensional
component of ${\hM }(x^{'}, y^{'})$ is compact up to the splitting of two
isolated trajectories for $J \in {\J }_*(L, \p)$.
\end{pro} \label{compa}

Proposition~\ref{bub} plays the key role in showing that $\d \circ \d = 0$ in 
Theorem~\ref{focohomology}.  We follow Floer-Oh's argument to show that 
${\bd }^{(r)} \circ {\bd }^{(r)} = 0$. The condition $\Si (L) \geq 3$, rather
than $\Si (L) \geq 2$, enters only in proving that $< \d \circ \d x, x> = 0$.
For $\Si (L) = 2$, Oh evaluates a number $\pmod 2$ of $J$-holomorphic
discs with Maslov index $2$ that pass through 
$x \in L \subset P$ in \cite{oh4}. 
Oh in \cite{oh4} verified that the number is always even, hence 
$< \d \circ \d x, x> = 0 \pmod 2$. In our case, this reflects to understand
the two lifts $x^{(r)}, g(x^{(r)})$ of $x$ with $ deg (g) = \pm 1$.
Note that $x^{(r)} \in (r , r + \s (L))$, $g(x^{(r)}) \in 
(r + deg (g) \s (L), r+ (deg (g) +1) \s (L))$. The integer graded symplectic
coboundary is not well-defined in this situation. We leave it to a future study.

\begin{lm} Under the same hypothesis in Proposition~\ref{bub}, 
$\bd^{(r)} \circ \bd^{(r)} = 0$.
\end{lm}
Proof: If $x \in C^{(r)}_{n-1}(L, \p; P, J)$, then by definition 
the coefficient of $z \in C^{(r)}_{n+1}(L, \p; P, J)$  in 
$\bd^{(r)} \circ \bd^{(r)}(x)$ is 
\begin{equation}
 \sum_{y \in C^{(r)}_{n}(L, \p; P, J)} \# {\hM }(x, y) \cdot \# {\hM }(y,z). 
\end{equation}
By Proposition~\ref{bub}, the boundary of the 1-dimensional
manifold $\hM (x, z) = {\CM }(x, z)/{\R}$ corresponds to two isolated
trajectories
${\CM }(x, y) \x {\CM }(y, z)$. Each term $\# \hM (x, y) \cdot \# \hM(y,z)$ is
the number of the 2-cusp trajectory of $\hM (x, z)$ with 
$y \in C^{(r)}_{n}(L, \p; P, J)$. 
For any such $y$ there are $J$-holomorphic curves $u \in {\CM }(x, y)$ and
$v \in {\CM }(y,z)$.
The other end of the corresponding 
component of  the 1-manifold $\hM (x,z)$ corresponding to the splitting
${\CM }(x, y^{'}) \x {\CM }(y{'},z)$ with $u^{'} \in {\CM }(x, y^{'})$ and 
$v^{'} \in {\CM }(y{'},z)$. Then $\hM (x, z)$ has an one
parameter family of paths 
from $x$ to $z$ with ends $u \# v$ and $u^{'} \# v^{'}$ for appropriate
grafting \cite{fl2}. If we lift $u$ to $\tilde{u} \in \tilde{\CM} (x^{(r)}, 
\tilde{y})$, then
\begin{equation} \label{rest}
1 = \mu_{\tilde{u}} = \mu^{(r)}(\tilde{y}) - \mu^{(r)}(x) = 
\mu^{(r)}(\tilde{y}) - (n-1).
\end{equation}
So $\mu^{(r)}(\tilde{y}) = n$; and $\tilde{y} = y^{(r)}$ is the preferred
lift, thus we have $\tilde{u} \in \tilde{\CM}(x^{(r)}, 
y^{(r)})$. Similarly $\tilde{v} \in \tilde{\CM}(y^{(r)}, z^{(r)})$.
Since $u^{'} \# v^{'}$ is homotopic to $u \# v$ rel $(x^{(r)}, z^{(r)})$, the 
lift $\tilde{u}^{'} \# \tilde{v}^{'}$ is also a path with ends 
$(x^{(r)}, z^{(r)})$. Now using the symplectic action $a$ is 
non-increasing along the gradient trajectory $\tilde{u}^{'}$, we have
\begin{equation}
r < a(z^{(r)}) \leq a(\tilde{y}^{'}) \leq a(x^{(r)}) < r + \s (L).
\end{equation} \label{ineq}
By uniqueness, $\tilde{y}^{'} = (y^{'})^{(r)}$ and using (\ref{rest})
for $u^{'}$, we have $\mu^{(r)}((y^{'})^{(r)}) = \mu^{(r)}(x^{(r)}) +1 = n$;
so $(y^{'})^{(r)} \in C^{(r)}_{n}(L, \p; P, J)$. Thus the number of
two-trajectories connecting $x^{(r)}$ and $z^{(r)}$ with index 2 is even. 
\qed

Now $(C^{(r)}_{n}(L, \p; P, J), \bd^{(r)}_n)_{n\in \Z}$ is indeed an integer
graded Floer cochain complex. We call its cohomology to be an integer
graded symplectic Floer cohomology,
denote by
\[I^{(r)}_*(L, \p; P, J) = H^*(C^{(r)}_*(L, \p; P, J), \bd^{(r)}), \ \ \
* \in \Z. \]
From the construction we see that if $[r,s] \subset {\R }_{L, \p}$, then
$I^{(r)}_*(L, \p; P, J) = I^{(s)}_*(L, \p; P, J)$.
The relation between $I^{(r)}_*(L, \p; P, J)$ and $HF^*(L, \p; P)$ will be 
discussed in \S 4.

\section{Invariance property of the integral symplectic Floer cohomology}

In this section we are going to show that the integer graded symplectic Floer 
cohomology defined in the previous section is invariant under the change of
$J$ and under the exact deformations of Lagrangian sub-manifolds.
\vspace{0.1in}

We consider an one parameter family $\{(J^{\lam }, \p^{\lam }\}_{\lam \in \R}$
that interpolates from $(J^0, \p^0)$ to $(J^1, \p^1)$ and is constant in
$\lam $ outside $[0,1]$. Let ${\P }_{1, \ve /2}$ 
be the set of all one parameter families
where $J^0, J^1 \in {\J }_{*}(L, \p) \subset {\J }_{reg}(L, \p)$ such that
\begin{equation}
  \int_{{\T } \times [0,1]}  | (J^{\lam}_{\t} - J^i) \frac{\bd u_{\lam }}
{\bd t} |^2 dt d\t  < \frac{\ve }{2} , \ \ \ 
|\p^{\lam } - id|_{C^1(\T \x [0,1])} < \frac{\ve }{2} ,
\end{equation}
for all $u_{\lam } \in {\P }(L, \p_{\lam }; P)$, $i = 0, 1$. Here we also assume
that $\p^{\lam }$ is  exact under the change of $\lam $, so $J^{\lam }_t 
= J(\lam , t), \p^{\lam }_t = \p (\lam, t), \p (\lam , 0) = Id$ are two
parameter families of almost complex structures compatible to $\o $ and
exact isotopes contractible to the identity. Such $\p^{\lam }_t$ connecting 
$\p^0, \p^1$ does exist.
Floer discussed the invariance of the Floer cohomology
under the change of $(J, \p)$ for
$(J^0, \p^0)$ $C^{\infty }$-close to $(J^1, \p^1)$ in \cite{fl3}.
 For $(J, \p )$, define
\[ C^{(r)}_{J, L} = \min \{ a(x^{(r)}) -r, \s (L) +r - a(x^{(r)}) | x \in
Z_{\p } \}.  \]
Since $r \in {\R }_{L, \p }$ is an regular value of the symplectic action and
modulo $\s (L) \Z$ the image $Im a (Z_{\p })$ is a finite set, so we have
$C^{(r)}_{J, L}$ is a positive number in $(0, \s (L))$. For instance, we may
choose $\ve = C^{(r)}_{J, L}/16$. 
\medskip

Define the perturbed $J$-holomorphic curve equations
\begin{equation} \label{jl}
{\pj }_{\lam } u_{\lam } (\t, t ) = 
\frac{\bd u_{\lam }}{\bd \t} + J^{\lam}_t 
\frac{\bd u_{\lam }}{\bd t} = 0 ,
\end{equation}  
with the moving Lagrangian coboundary conditions
\begin{equation} \label{ll}
u_{\lam }(\t, 0) \in L,  u_{\lam }(\t, 1) \in \p^{\lam }_1(L) . 
\end{equation} 
This directly generalizes the $J$-holomorphic curve equation in the case
of $(J^0, \p^0)$ and $(J^1, \p^1)$. The moduli space ${\CM }_{\lam } (x, y )$
of (\ref{jl}) and (\ref{ll}) has the same analytic properties as the 
moduli space ${\CM }(x, y)$ except for the translational invariance (see
Proposition 3.2 in \cite{fl2}). Hofer in \cite{ho} anaylized the compactness 
property for a similar moving Lagrangian coboundary condition, Oh in \cite{oh}
studied that the bubbling-off a sphere or disk can not occur in the components
of ${\CM }_{\lam } (x, y )$ for the monotone Lagrangian sub-manifold $L$ and
${\Si }(L) \geq 3$. The index of $u_{\lam }$ can be proved to be the same 
as a topological index for the moduli space of perturbed $J$-holomorphic curves.
The arguments of the proof of invariance under the change of $(J, \p)$ are 
the same as in \cite{fl2}, \cite{fl3} and \cite{oh}. But we need to ensure that
the cochain map is well-defined for the integer graded cochain complexes.
\begin{lm} \label{lift}
 If $u_{\lam} \in {\CM }_{\lam} (x_0, x_1)$, 
$(J^{\lam }_t, \p^{\lam }_t) \in {\P }_{1, \ve /2}$ and 
$\tilde{u}_{\lam } \in {\P }(\tilde{x}_0, \tilde{x}_1)$ 
is any lift of $u_{\lam }$, then
\[a_{(J^1, \p^1)}(\tilde{x}_1) < a_{(J^0, \p^0)} (\tilde{x}_0) + \ve . \]
\end{lm} 
Proof: Note that the path $\{u_{\lam }(\t ) | \t \in (- \infty , 0) \}$ is
a gradient trajectory for $(J^0, \p^0)$ and $\{u_{\lam }(\t ) | \t \in 
(1, \infty )\}$ is a gradient trajectory for $(J^1, \p^1)$. So we have
\begin{equation} \label{ine}
 a_{(J^0, \p^0)} (\tilde{u}(0)) \leq a_{(J^0, \p^0)} (\tilde{x}_0), 
\ \ \ a_{(J^1, \p^1)}(\tilde{x}_1) \leq a_{(J^1, \p^1)}(\tilde{u}(1)) .
\end{equation} 
Since $u_{\lam } \in {\CM }_{\lam }(x_0, x_1)$, by the property of 
${\P }_{1, \ve /2}$, we have
\begin{eqnarray*}
I_{\o }(u_{\lam})|_{\T \x [0,1]} & = & 
\|\bd_{J^{\lam }} u_{\lam }\|^2_{L^2(\T \x [0,1])} - 
\|{\pj }_{\lam }u_{\lam}\|^2_{L^2(\T \x [0,1])} \\ 
& \geq & \|\bd_{J^0} u_{\lam }\|^2_{L^2(\T \x [0,1])} -
\|(\bd_{J^{\lam }} - \bd_{J^0})u_{\lam }\|^2_{L^2(\T \x [0,1])}\\
 & \geq & - \frac{1}{2}
\|(J^{\lam }_t - J^0) \frac{\bd u_{\lam }}{\bd t} \|^2_{L^2(\T \x [0,1])}\\
 & > & - \frac{\ve }{2}
\end{eqnarray*}
Note that 
$\bd_J u = \frac{1}{2} ( \frac{\bd u}{\bd \t} - J \frac{\bd u}{\bd t})$. 
So the symplectic actions for $\tilde{u}(0)$ and $\tilde{x}_1$ are related
by the following.
\begin{equation} \label{gei}
 a(\tilde{u}(0)) - a(\tilde{x}_1) = I_{\o }(u_{\lam})|_{\T \x [0,1]} > 
- \frac{\ve }{2} . \end{equation}
Thus from (\ref{ine}) and (\ref{gei}) we get 
\begin{equation} \label{oneway}
a_{(J^0, \p^0)} (\tilde{x}_0) \geq a_{(J^0, \p^0)} (\tilde{u}(0)) > 
a(\tilde{x}_1) - \frac{\ve }{2} .
\end{equation} 
Then the result follows. \qed

\begin{df} \label{cochain}
For each $n$, define a homomorphism
$\p_{01}^n : C^{(r)}_n (L, \p^0; P, J^0) \to C^{(r)}_n (L, \p^1; P, J^1)$ by
\[ \p_{01}^n (x_0) = \sum_{x_1 \in C^{(r)}_n (L, \p^1; P, J^1)} \#
{\CM }^0_{\lam } (x_0, x_1) \cdot x_1. \]
\end{df}

We show that the homomorphism $\p_{01}^n$ is indeed a cochain map with respect
to the integral lift $r$.

\begin{lm} \label{01map}
$\{\p_{01}^*\}_{* \in Z}$ is a cochain map, i.e. 
\[ \bd^{(r)}_{n,1} \circ \p_{01}^n = \p_{01}^n \circ \bd^{(r)}_{n,0},\]
for all $n \in Z$. 
\end{lm}
Proof: For $x_0 \in C^{(r)}_n (L, \p^0; P, J^0)$ and 
$y_1 \in C^{(r)}_{n+1} (L, \p^1; P, J^1)$,the coefficient of $y_1$ in
$(\bd^{(r)}_{n,1} \circ \p_{01}^n  - \p_{01}^n \circ \bd^{(r)}_{n,0})(x_0)$
is the number of the set modulo 2,
\begin{equation} \label{42}
 \cup_{y_0 \in C^{(r)}_{n+1} (L, \p^0; P, J^0)} - \hat{\CM }_{J^0, \p^0}
(x_0, y_0) \x {\CM }_{\lam }^0 (y_0, y_1) \cup_{x_1 \in
C^{(r)}_n (L, \p^1; P, J^1)} {\CM }_{\lam }^0 (x_0, x_1) \x
\hat{\CM }_{J^1, \p^1} (x_1, y_1) .
\end{equation} 
The ends of the 1-dimensional manifold
${\CM }_{\lam }^1(x_0, y_1)$ are in one-to-one correspondence with the set
\begin{equation} \label{43}
 \cup_{y \in Z_{{\p }_0 }} - {\hM }_{J^0, \p^0}
(x_0, y) \x {\CM }_{\lam }^0 (y, y_1) \cup_{x \in Z_{{\p }_1}}
{\CM }_{\lam }^0 (x_0, x) \x 
{\hM }_{J^1, \p^1} (x, y_1) . \end{equation}
For an end $u \# v$ of ${\CM }_{\lam }^1(x_0, y_1)$ corresponding to an 
element in (\ref{42}), the other end of the component $u^{'} \# v^{'}$ 
corresponds to an element in (\ref{43}). For $u \in {\CM }^0_{\lam }(x_0, y)$
and $v \in {\hM }^1_{J^1, \p^1}(y, y_1)$, ${\CM }^1_{\lam }(x_0, y_1)$
gives a 1-parameter family of paths in ${\P }(L, \p_{\lam }; P)$ with 
fixed end points $x_0, y_1$. Such a path gives a homotopy of paths from
$u \#_{\r } v$ to $u^{'} \#_{\r } v^{'}$ rel end points. The lift of 
$u \#_{\r } v$ starting $x_0^{(r)}$ ending at $y_1^{(r)}$, so the same
is true for the lift of $u^{'} \#_{\r } v^{'}$ starting at $x_0^{(r)}$. 
Suppose $u^{'}$ lifts to an element in ${\CM }^0_{\lam }(x_0^{(r)}, \tilde{y})$.
By Lemma~\ref{lift}, we have 
\[ a_{(J^1, \p^1)} (\tilde{y}) < a_{(J^0, \p^0)} (x_0^{(r)}) + \ve 
< r + \s (L). \] 
Furthermore, 
using the trajectory decreasing the symplectic action, we also have
\[ a_{(J^1, \p^1)} (\tilde{y}) > a_{(J^1, \p^1)} (y_1^{(r)}) > r . \]
So this shows that $a_{(J^1, \p^1)} (\tilde{y}) \in (r, r+\s (L))$ 
which gives the
preferred lift $\tilde{y} = y^{(r)}$. Now using Proposition~\ref{mas}, 
\[ 1 = \mu^{(r)}(y_1^{(r)}) - \mu^{(r)}(y^{(r)}) = (n+1) - \mu^{(r)}(y^{(r)}),\]
so we get $\mu^{(r)}(y^{(r)}) = n$, $y \in C_n^{(r)}(L, \p^1; P, J^1)$. This
shows that $u^{'} \#_{\r } v^{'}$ in (\ref{43}) actually corresponds to an 
element in (\ref{42}). So the cardinality is always even. \qed

For $(J^i, \p^i) \in {\P }_1$, $i = 0, 1, 2$, and define 
new classes of perturbations ${\P }_{2, \ve}$ to consists of
\[ (J^{\lam }, {\p }^{\lam }) = \left\{ 
\begin{array}{ll}
(J^0, \p^0) & \lam \leq - T \\
(J^1, \p^1) & - T +1 \leq \lam \leq T -1 \\
(J^2, \p^2) & \lam \geq T
\end{array} 
\right. \]
for some fixed positive number $T ( > 2)$, such that
\begin{equation} \label{p2}
\int_{{\T } \x ([-T, -T +1] \cup [T-1, T])} |(J^{\lam } - J^i)
\frac{ \bd u_{\lam }}{\bd t} |^2 dt d \t < \ve .
\end{equation}

Then if $(J^{\lam }_1, {\p }^{\lam }_1) \in 
{\P }_{1, \ve /2}((J^0, \p^0), (J^1, \p^1))$
and $(J^{\lam }_2, {\p }^{\lam }_2) \in 
{\P }_{1, \ve /2}((J^1, \p^1), (J^2, \p^2))$,
then we can compose $(J^{\lam }_1, {\p }^{\lam }_1)$ with 
$(J^{\lam }_2, {\p }^{\lam }_2)$ to get $(J^{\lam }, {\p }^{\lam }) \in
{\P }_{2, \ve }((J^0, \p^0),(J^2, \p^2))$. We denote such a composition
as $(J^{\lam }, {\p }^{\lam }) =
(J^{\lam }_1, {\p }^{\lam }_1) \#_T (J^{\lam }_2, {\p }^{\lam }_2)$.
Then for a large fixed $T$ and each compact set $K$ in
${\CM }_{(J^{\lam }_1, {\p }^{\lam }_1)}(x, y) \x
{\CM }_{(J^{\lam }_2, {\p }^{\lam }_2)}(y, z)$, there is a ${\r }_T > 0$ and 
for all $\r > {\r }_T$ a local diffeomorphism
\[ \# : {\CM }_{(J^{\lam }_1, {\p }^{\lam }_1)}(x, y) \x
{\CM }_{(J^{\lam }_2, {\p }^{\lam }_2)}(y, z) \supset K \to 
{\CM }_{(J^{\lam }_1, {\p }^{\lam }_1) \#_T (J^{\lam }_2, {\p }^{\lam }_2)}
(x, z) . \]
See Proposition 2d.1 in \cite{fl3}. 
\begin{lm} \label{compose}
For the above $(J^{\lam }, {\p }^{\lam })$ and $\r > {\r }_T$, 
\[ \p_{02}^n = \p_{12}^n \circ \p_{01}^n , \]
for all $n \in \Z$.
\end{lm}
Proof: For $x_0 \in C_n^{(r)}(L, \p^0; P, J^0)$, 
\[ \p_{02}^n (x_0 ) = \sum_{y_0 \in C_n^{(r)}(L, \p^2; P, \J^2)} \# 
{\CM }^0_{\lam }(x_0, y_0) \cdot y_0, \]
\[ \p_{12}^n \circ \p_{01}^n (x_0 ) = \sum \# 
({\CM }^0_{(J^{\lam }_1, {\p }^{\lam }_1)}(x_0, y) \x 
{\CM }_{(J^{\lam }_2, {\p }^{\lam }_2)}(y,y_0) ) \cdot y_0, \]
where the summation $\sum $ runs over $y \in C_n^{(r)}(L, \p^1; P, J^1)$
and $y_0 \in C_n^{(r)}(L, \p^2; P, J^2)$. The local diffeomorphism $\#$ 
implies that 
\[ \# {\CM }^0_{\lam }(x_0, y_0) = 
\# ({\CM }^0_{(J^{\lam }_1, {\p }^{\lam }_1)}(x_0, y) \x 
{\CM }^0_{(J^{\lam }_2, {\p }^{\lam }_2)}(y,y_0)) .\]
All we need to check is that $y \in C_n^{(r)}(L, \p^1; P, J^1)$, 
this can be checked by 
the same argument in Lemma~\ref{01map}. \qed

For two classes $(J^{\lam }, {\p }^{\lam })$, $(\overline{J}^{\lam }, 
\overline{{\p }}^{\lam })$ in
${\P }_{2, \ve}((J^0, \p^0),(J^2, \p^2))$, the following lemma shows that the 
induced cochain maps $\p_{02}^n$, $\overline{\p }_{02}^n$ are cochain homotopic
to each other.

\begin{lm} \label{homotopy}
If $(J^{\lam }, {\p }^{\lam })$, $(\overline{J}^{\lam }, 
\overline{{\p }}^{\lam })$ in
${\P }_{2, \ve} ((J^0, \p^0),(J^2, \p^2))$ 
can be smoothly deformed from one to another
by a 1-parameter family $(J^{\lam }_s, {\p }^{\lam }_s), s \in [0, 1]$, i.e.
$(J^{\lam }_s, {\p }^{\lam }_s) = (J^{\lam }, {\p }^{\lam })$ for
$s \leq 0$ and $(J^{\lam }_s, {\p }^{\lam }_s) =  (\overline{J}^{\lam }, 
\overline{{\p }}^{\lam })$ for $s \geq 1$. Then $\p_{02}^*$, 
$\overline{\p }_{02}^*$ are cochain homotopic to each other.
\end{lm}
Proof: we need to construct a homomorphism
\[ H: C^{(r)}_*(L, \p^0; P, J^0) \to C^{(r)}_*(L, \p^2; P, J^2) , \]
of degree $- 1$ with the property
\begin{equation} \label{hom}
 \p_{02}^n - \overline{\p }_{02}^n = H {\bd }_{n,0}^{(r)} + 
{\bd }_{n,2}^{(r)} H, \end{equation}
for all $n \in \Z$. Associated to $(J^{\lam }_s, {\p }^{\lam }_s)$, there is 
a moduli space $H{\CM }(x_0, y_0) = \cup_{s \in [0,1]}
{\CM }^0_{(J^{\lam }_s, {\p }^{\lam }_s)}(x_0, y_0)$.
\[ H{\CM }(x_0, y_0) = \{ (u, s) \in 
{\CM }^0_{(J^{\lam }_s, {\p }^{\lam }_s)}(x_0, y_0) \x [0, 1] \} \subset 
{\P }(L, {\p }^{\lam }_s; P)(x_0, y_0) \x [0,1]. \]
$H{\CM }(x_0, y_0)$ are regular zero sets of 
${\pj }_{(J^{\lam }_s, {\p }^{\lam }_s)}$ and are smooth manifolds of
dimension $\mu^{(r)}(y_0^{(r)}) - \mu^{(r)}(x_0^{(r)}) + 1$.
Consider the case of $\mu^{(r)}(x_0^{(r)}) = \mu^{(r)}(y_0^{(r)}) = n$, 
the boundary of 1-dimensional sub-manifold $H{\CM }(x_0, y_0)$ of 
${\P }(L, {\p }^{\lam }_s; P)(x_0, y_0) \x [0,1]$ consists of
\begin{itemize}
\item ${\CM }^0_{(J^{\lam }, {\p }^{\lam })}(x_0, y_0) \x \{0 \} 
\cup {\CM }^0_{(\overline{J}^{\lam }, 
\overline{{\p }}^{\lam })}(x_0, y_0) \x \{1 \} $
\item $\cup_{s\in [0, 1], y} 
{\CM }^{-1}_{(J^{\lam }_s, {\p }^{\lam }_s)}(x_0, y)
\x {\CM }^0_{(J^2, \p^2)}(y, y_0)$ for $y \in C_{n-1}^{(r)}(L, \p^2; P, J^2)$.
\item ${\CM }^0_{(J^0, \p^0)}(x_0,x) \x \cup_{s\in [0, 1], x}
{\CM }^{-1}_{(J^{\lam }_s, {\p }^{\lam }_s)}(x, y)$ for
$x \in C_{n-1}^{(r)}(L, \p^0; P, J^0)$.
\end{itemize}

Note that ${\CM }^{-1}_{(J^{\lam }_s, {\p }^{\lam }_s)}(x_0, y)$ and 
${\CM }^{-1}_{(J^{\lam }_s, {\p }^{\lam }_s)}(x, y)$ are solutions of
$(u, s)$ of $J$-holomorphic equations lying in virtual dimension $-1$, they can
only occur  for $0 < s < 1$. Define $H: C^{(r)}_n(L, \p^0; P, J^0) \to 
C^{(r)}_{n-1}(L, \p^2; P, J^2)$ by
\begin{equation}
H(x_0) = \sum_y \sum_s \# 
{\CM }^{-1}_{(J^{\lam }_s, {\p }^{\lam }_s)}(x_0, y).
\end{equation} 
Similar to Lemma~\ref{01map}, by checking the corresponding preferred lifts and
the integral Maslov index, we get the desired cochain homotopy $H$ between
${\p }_{02}^n$ and $\overline{\p }_{02}^n$ satisfying (\ref{hom}). \qed

From Lemma~\ref{homotopy}, we have that $\p_{02}^*$ from $(J^{\lam }, 
{\p }^{\lam })$ is the same homomorphism of $\overline{\p }_{02}^*$ from
$(\overline{J}^{\lam }, \overline{\p }^{\lam })$ on the integer
graded symplectic Floer cohomology. Then we can prove the invariance of the
integer graded symplectic Floer cohomology under the continuation of $(J,\p)$.

\begin{thm} \label{invariant}
For any continuation $(J^{\lam }, {\p }^{\lam }) \in {\P }_{1, \ve /2}$ which 
is regular at the ends, there exists an isomorphism
\[\p_{02}^n: I^{(r)}_n(L, \p^0; P, J^0) \to I^{(r)}_n(L, \p^1; P, J^1), \]
for all $n \in \Z$.
\end{thm}
Proof: Let $(J^{- \lam }, {\p }^{ - \lam })$ be the reversed family 
of $(J^{\lam }, {\p }^{\lam })$ by $\t = - {\t }^{'}$. So we can form
a family of composition 
$(J^{\lam }, {\p }^{\lam }) \#_T (J^{- \lam }, {\p }^{ - \lam })$ in 
${\P }_{2, \ve}$ for some fixed $T ( > 2)$. By Lemma~\ref{compose}, we have
$${\p }_{(J^{\lam }, {\p }^{\lam }) \#_T (J^{- \lam }, {\p }^{ - \lam })} =
{\p }_{10}^* \circ {\p }_{01}^* . $$
For $(J^{\lam }, {\p }^{\lam }) \#_T (J^{- \lam }, {\p }^{ - \lam })$, it can be
deformed to the trivial continuation as $(J^0, \p^0)$ for all $\t \in \R$.
Then by Lemma~\ref{homotopy}, we have 
\[ {\p }_{10}^* \circ {\p }_{01}^* = {\p }_{00}^* = id:
I^{(r)}_*(L, \p^0; P, J^0) \to I^{(r)}_*(L, \p^0; P, J^0). \]
Similarly, ${\p }_{01}^* \circ {\p }_{10}^* = {\p }_{11}^* = id$ on
$I^{(r)}_*(L, \p^1; P, J^1)$. Thus the result follows. \qed

Now the integer graded symplectic Floer cohomology $I_*^{(r)}$ is functorial 
with respect to compositions of continuations $(J^{\lam }, {\p }^{\lam })$
and invariant under continuous deformations of $(J^{\lam }, {\p }^{\lam })$
within the set of continuations ${\P }_{1, \ve /2}$.

\section{Spectral sequence for the symplectic Floer cohomology}

In this section we are going to show that the $Z$-graded symplectic Floer
cohomology $I_*^{(r)}(L, \p;P)$ for $r \in {\R }_{L, \p}$ 
and $* \in \Z$ determines
the ${\Z }_{ \Si (L) }$-graded symplectic Floer cohomology 
$HF^*(L, \p; P), * \in 
{\Z }_{\Si (L)}$. The way to link them together is to filter the integer
graded Floer cochain complex. The filtration, by a standard method,
 formulates a spectral sequence
which converges to the ${\Z }_{\Si (L)}$ graded symplectic Floer cohomology 
$HF^*(L, \p; P)$.
\vspace{0.1in}

\begin{df} \label{filter}
For $r \in {\R }_{L, \p}, j \in {\Z }_{\Si (L) }$ and 
$n \equiv j \ \ (\ mod \ \ \Si (L))$, define the free $Z_2$ modules
\[ F_n^{(r)} C_j(L, \p; P, J) = \sum_{k \geq 0} C^{(r)}_{n + \Si (L) k}
(L, \p; P, J), \]
which gives a natural decreasing filtration on $C_*(L, \p; P, J), * \in 
{\Z }_{\Si (L) }$. 
\end{df}

There is a finite length decreasing filtration of $C_j(L, \p; P, J)$, 
$j \in {\Z }_{\Si (L) }$:
\begin{equation} \label{dea}
\cdots \subset F_{n+\Si (L)}^{(r)}C_j(L, \p; P, J) \subset F_n^{(r)}
  C_j(L, \p; P, J) \subset F_{n - \Si (L) }^{(r)}C_j(L, \p; P, J)
  \subset  \cdots \subset C_j(L, \p; P, J) .
\end{equation}
\begin{equation} \label{col}
C_j(L, \p; P, J) = \cup_{n \equiv j (mod \ \ \Si (L))} 
F_n^{(r)} C_j(L, \p; P, J).
\end{equation}
Note that the symplectic action is non-increasing 
along the gradient trajectories
($J$-holomorphic curves on $\T \x \R$), it follows that the coboundary map
$\d: F_n^{(r)} C_j(L, \p; P, J) \to F_{n+1}^{(r)} C_{j+1}(L, \p; P, J)$ 
(in Theorem~\ref{focohomology})
preserves the
filtration in definition~\ref{filter}. Thus the ${\Z }_{\Si (L)}$
graded symplectic Floer cochain complex $(C_j(L, \p; P, J), \d)_{j \in 
{\Z }_{\Si (L)}}$ has a decreasing bounded 
filtration $(F_n^{(r)}C_*(L, \p; P, J), \d)$,
\begin{equation} \label{pre}
\begin{array}{ccccc}
 & \downarrow & \downarrow &  & \downarrow \\
\cdots & \subset F_{n+\Si (L)}^{(r)}C_j(L, \p; P, J) & 
\subset F_n^{(r)}C_j(L, \p; P, J) & \subset
\cdots & \subset C_j(L, \p; P, J) \\
 & \downarrow \partial^{(r)} & \downarrow \partial^{(r)} &  & \downarrow \d \\
\cdots & \subset F_{n+\Si (L) +1}^{(r)}C_{j+1}(L, \p; P, J) & 
\subset F_{n+1}^{(r)}C_{j+1}(L, \p; P, J)
& \subset \cdots & \subset C_{j+1}(L, \p; P, J) \\
& \downarrow & \downarrow &  & \downarrow
\end{array} .
\end{equation}
\begin{lm} (1) The cohomology of the vertical cochain subcomplex 
$F_n^{(r)} C_*(L, \p; P, J)$ in the filtration (\ref{pre}) is 
$F_n^{(r)} I_j^{(r)}(L, \p; P, J).$

(2) There is a natural bounded filtration for  
$\{I_*^{(r)}(L, \p; P, J) \}_{* \in Z}$ the integer graded symplectic
Floer cohomology,
\[\cdots F_{n+\Si (L)}^{(r)}HF^j(L, \p; P, J) \subset F_n^{(r)} 
HF^j(L, \p; P, J) \subset F_{n-\Si (L)}^{(r)}HF^j(L, \p; P, J) 
\cdots \subset I_j^{(r)}(L, \p; P, J),\] 
where $F_n^{(r)}HF^j(L, \p; P, J) = \ker (I_j^{(r)}(L, \p; P, J)) \to
F_{n- \Si (L)}^{(r)}I_j^{(r)}(L, \p;P, J))$.
\end{lm}
Proof: The results follows from definition~\ref{filter} and standard results
in \cite{sp} Chapter 9. \qed

\begin{thm} \label{E1} For $\Si (L) \geq 3, $

(i) There is a spectral sequence $(E^k_{n, j}, d^k)$ with
\[ E^1_{n,j}(L, \p; P, J) \cong I_n^{(r)}(L, \p; P, J), \ \ \ 
n\equiv j \ \ (mod \ \ \Si (L)), \]
and 
\[E^{\infty }_{n,j}(L, \p; P, J) \cong F_n^{(r)} HF^j(L, \p; P, J)/ 
F_{n + \Si (L)}^{(r)}HF^j(L, \p; P, J) . \]

(ii) The spectral sequence $(E^k_{n, j}, d^k)$ converges to the $Z_{\Si (L)}$
graded symplectic Floer cohomology $HF_*(L, \p; P, J)$, where
\[ d^k: E^k_{n, j}(L, \p; P, J) \to 
E^k_{n + \Si (L) k +1, j +1}(L, \p; P, J) . \] 
\end{thm}
Proof: (i) Note that 
\[ F_n^{(r)} C_j(L, \p; P, J) /F_{n + \Si (L)}^{(r)} C_j(L, \p; P, J) = 
C_n ^{(r)}(L, \p; P, J) . \]
It is standard from \cite{sp} that there exists a spectral sequence 
$(E^k_{n, j}, d^k)$ with $E^1$ term given by the cohomology of $
F_n^{(r)} C_j(L, \p; P, J) /F_{n + \Si (L)}^{(r)} C_j(L, \p; P, J)$, so we have 
$E^1_{n,j}(L, \p; P, J) \cong I_n^{(r)}(L, \p; P, J)$
and $E^{\infty }_{n,j}(L, \p; P, J)$ is isomorphic to the bigraded $Z_2$-module
associated to the filtration $F^{(r)}$ of the $Z$-graded symplectic
Floer cohomology $I_n^{(r)}(L, \p; P, J)$.

(ii) Since the Lagrangian intersections are transversal and $Z_{\p}$ is a finite
set of intersections $L \cap \p (L)$, so the filtration $F$ is bounded
and complete from (\ref{col}). Thus the spectral sequence converges to the
$Z_{\Si (L)}$ graded symplectic Floer cohomology. 
Note that the grading is unusual
(jumping by $\Si (L)$ in $n$), we list the terms for $Z^k_{*, *}, E^k_{*,*}$.
\[ Z^k_{n,j}(L, \p; P, J) = \{ x \in F_n^{(r)}C_j(L, \p; P, J) | \d x \in
F_{n+1 + \Si (L) k}^{(r)} C_{j+1}(L, \p; P, J) \} \]
\[E^k_{n,j}(L, \p; P, J) = Z^k_{n,j}(L, \p; P, J)/
\{ Z^{k+1}_{n + \Si (L),j}(L, \p; P, J) + 
\d Z^{k-1}_{n+ (k-1) \Si (L) -1,j-1}(L, \p; P, J)\} \]
\[Z^{\infty}_{n,j}(L, \p; P, J) = 
\{ x \in F_n^{(r)}C_j(L, \p; P, J)  | \d x = 0 \} \]
\[E^{\infty}_{n,j}(L, \p; P, J) = Z^{\infty}_{n,j}(L, \p; P, J) /
\{ Z^{\infty}_{n + \Si (L),j}(L, \p; P, J) +
d Z^{\infty}_{n+ (k-1) \Si (L) -1,j-1}(L, \p; P, J)\} . \]
Thus $\d$ induces the higher differential
\[ d^k: E^k_{n, j}(L, \p; P, J) \to 
E^k_{n + \Si (L) k +1, j +1}(L, \p; P, J) . \] \qed

\begin{thm} For $\Si (L) \geq 3$, 

(1) for any continuation $(J^{\lam }, {\p }^{\lam }) \in {\P }_{1, \ve /2}$ 
which is
regular at ends, there exists an isomorphism
\[E^1_{n,j}(L, {\p}^0; P,J^0) \cong E^1_{n,j}(L, {\p}^1; P,J^1) .\]

(2) For each $k \geq 1$, $E^k_{n, j}(L, \p; P, J)$ are the symplectic invariant
under continuous deformations of $(J^{\lam }, {\p }^{\lam })$ within the
set of continuations.
\end{thm}
Proof: Clearly (2) follows from (1) by Theorem 1 in \cite{sp} page 468.
From Theorem~\ref{E1}, we have an isomorphism
$E^1_{n,j}(L, {\p}^0; P,J^0) \cong I_n^{(r)}(L, {\p}^0; P,J^0)$, 
so there exists an isomorphism from \S 3, 
$I_n^{(r)}(L, {\p}^0; P
,J^0) \to I_n^{(r)}(L, {\p}^1; P,J^1)$ which respects the filtration and 
induces an isomorphism on the $E^1$ term. \qed

Now we can view that $E^k_{n, j}(L, \p; P, J) = E^k_{n, j}(L, \p; P)$, 
$E^1_{n,j}(L,\p; P) = I_n^{(r)}(L,\p; P)$, for $k \geq 1, r \in {\R }_{L, \p}$
, are new symplectic invariants provided $\Si (L) \geq 3$.
All these new symplectic invariants should contain more information on
$(P, \o; L, \p)$, they are also finer than $HF^*(L,\p; P),
 * \in {\Z }_{\Si (L)}$ the usual Floer cohomology. 
In particular, the minimal $k$ for which
$E^k = E^{\infty }$ should be meaningful, denoted by $k(L)$.
We will discuss the applications of these symplectic
invariants on Lagrangian embedding, Maslov index and Hofer symplectic 
energy norm. 

\begin{cor}
For $\Si (L) \geq 3, j \in {\Z }_{\Si (L)}$, 
\[ \sum_{k \in Z} I^{(r)}_{ j + \Si (L) k }(L, \p; P) = HF^j (L, \p; P) \]
if and only if all the differentials $d^k$ in the spectral sequence 
$(E^k_{n, j}, d^k)$ are trivial (i.e. $k(L) = 1$).
\end{cor}

In general, we see that $\sum_{k \in Z} I^{(r)}_{ j + \Si (L) k }(L, \p; P)
\neq HF^j (L, \p; P)$ for $j \in Z_{\Si (L)}$. $I_*^{(r)}(L,\p; P)$ 
can be thought as an integer lift of the 
symplectic Floer cohomology $HF^*(L, \p; P)$.

\section{Applications and Remarks}

\subsection{Hofer's energy and Chekanov's construction}

Hofer \cite{ho1} introduced the notion of disjunction energy or 
displacement energy associated with a subset of symplectic manifold. Roughly, 
the Hofer's symplectic energy measures how large a variation of a (compactly 
supported) Hamiltonian function must be in order to push the subset off itself
by a time-one map of corresponding Hamiltonian flow. 
Hofer showed that his symplectic energy of every
open subset in standard symplectic vector space is nontrivial.
For more geometric study of the Hofer's energy, we refer to \cite{lm1} and 
\cite{lm2}.

\begin{df} Let ${\cal H}$ be the space of compactly supported functions on
$[0,1] \x P$. The Hofer's symplectic energy of a symplectic diffeomorphism 
$\p : P \to P$ is defined by
\[E(\p ) = \inf \{ \int^1_0(\max_{x \in P} H(s, x) - \min_{x \in P} H(s, x))ds|
 \mbox{$\p $ is a time one flow generated by $H \in {\cal H}$} \}. \]
\[ e_H( L) = \inf \{ E( \p) : \p \in Ham (P),  L \cap \p (L) = \emptyset 
\ \ \ \mbox{empty set} \}.\]
\end{df}

\begin{thm}[Chekanov \cite{ch}] \label{hse}
If $E(\p ) < \s (L)$, $L$ is rational and
$L$ intersects $\p (L)$ transversally, then
\[ \# (L \cap \p (L) ) \geq SB(L; Z_2), \]
where $SB(L; Z_2)$ is the sum of Betti number with $Z_2$ coefficients. (I.e. 
$e_H(L) \geq \s (L)$.)
\end{thm}

\noindent{\bf Remarks:} (1) Polterovich \cite{po4} used Gromov's figure $8$ 
trick and a refinement of Gromov's existence scheme of $J$-holomorphic
disc to show that $e_H(L) \geq  \frac{\s (L)}{2}$. Chekanov extends the
result to $e_H(L) \geq \s (L)$ which is optimal for general Lagrangian
sub-manifolds, see \cite{ch} \S 1 and \cite{oh2} \S 6, and he asked whether 
Theorem~\ref{hse} remains true if $E(\p ) = \s (L)$. This is the case 
corresponding to $r_0 \equiv 0$ $(mod  \s (L))$ which is not an regular value
for the symplectic action $a$. 
(2) Sikorav \cite{si} showed that $e_H(T^m) \geq \s (T^m)$. Theorem~\ref{hse}
generalizes the Sikorav's result to all rational Lagrangian sub-manifolds.

In \cite{ch}, Chekanov was able to use a restricted Floer cohomology in the study
of Hofer's symplectic energy of rational Lagrangian sub-manifolds. Denote
\[ {\O }_s = \{ \g \in C^{\infty }([0,1], P) | \g (0) \in L, \g (1) \in \p_s (L)\}, \]
\[ \O = \cup_{s \in [0,1]} {\O }_s  \subset [0,1] \x C^{\infty }([0,1], P). \]
One may choose the anti-derivative of $Da_s(z) \xi$ as $a_s: {\O }_s \to 
{\R }/{\s (L) Z}$. Chekanov fixed $a_0$ with critical value $0$. 
Pick $z_s \in L \cap {\p}_s(L)$ such that modulo $\s (L) \Z$, 
\[0 < a_s (z_s) = \min \{a_s (x) (mod \s (L) {\Z}) | x \in Z_{\p_s} \} 
< \s (L).\]
This is possible since $a_s (x) \equiv 0 (mod \s (L) {\Z})$ for all $x \in
Z_{\p_s}$ will contradict with $E(\p ) < \s (L)$.
So we may take $r_0 \neq 0$ sufficiently small in 
${\R }_{L, {\p }_s} \cap (0, \s (L))$, say $0 < r_0 < \frac{1}{8} a_s(z_s)$. 
The condition $E(\p ) < \s (L)$ provides that 
there is a unique $x \in Z_{{\p }_s}$ (the $x^{(r_0)}$) which
correspond to the unique lifts in $(r_0, r_0 + \s (L))$, 
\begin{equation} 0 < a_s(x) - a_s(z_s) < \s (L). 
\end{equation}
Under these restrictions, define $C_s$ the
free $Z_2$ module generated by $L \cap \p_s (L)$, and a coboundary map
$\bd_s \in End(C_s)$ (see below) such that ${\bd }_s \circ {\bd }_s =0$, so
$H^*(C_s, {\bd }_s )$ is well-defined for every $s \in [0,1]$.
Note that with the unique lifting $x, y$ in $(r_0, r_0 + \s (L))$ we have
the Chekanov's restricted Floer homology in the integer grading.
\begin{lm} \label{cs}
For $r_0$ as above, we have $C^{(r_0)}_*(L, \p_s; P, J) = C_s$.
Let ${\CM }_s(x, y)$ be the restricted moduli space of $J$-holomorphic 
curves as the set $\{ u \in {\CM}(L, \p_s (L)) |
u^*(\o ) = a_s(x) - a_s(y) \}$. So we have
\[ \bd^{(r_0)} x = \bd_s x = \sum \# {\hM }_s(x, y) y . \]
\end{lm} 
Proof: Note that $\mu_u = \mu^{(r_0)}(y) - \mu^{(r_0)}(x) = dim {\CM}_s(x,y)$ 
for any $u \in {\CM}_s(x,y)$. So $y$ in the coboundary $\bd_s x$ is the element
in $C^{(r_0)}_{n+1}(L, \p_s; P,J)$; for any $u \in {\CM}(x,y)$ in the
coboundary of $\bd^{(r_0)} x$, we have $u^*(\o ) = a_s(x) - a_s(y)$ (\cite{fl2}
Proposition 2.3). 
For the unique lift in $(r_0, r_0 + \s (L))$, the choice of $a_s$ makes
that there is an one-to-one correspondence between ${\CM }_s(x, y)$ and
${\CM }(x, y)$ for $\mu^{(r_0)}(y^{(r_0)}) - \mu^{(r_0)}(x^{(r_0)}) = 1$. 
So the coboundary maps agree on the $Z_2$ coefficients. \qed

Hence for the choice of the $r_0$ we have the identification between
Chekanov's restricted cohomology $H^*(C_s, {\bd }_s)$ and
our integer graded Floer cohomology $I^{(r_0)}(L, {\p }_s; P)$.

For $s$ sufficiently small, $x \in L \cap {\p }_s(L)$ is also a critical
point of the Hamiltonian function $H_s$ of $\p_s$. Then the Maslov index is 
related to the usual Morse index of $H_s$ in the following:
\begin{equation}
\mu^{(r_0)}(x^{(r_0)}) = \mu_{H_s}(x) - m .
\end{equation}
See \cite{ar}, \cite{du}, \cite{hs} and \cite{sz}.

\begin{pro} \label{chek}
For $r_0$ as above, we assume that (i) $L$ is monotone Lagrangian sub-manifold
in $P$, (ii) $\Si (L) \geq 3$, (iii) $E(\p ) < \s (L)$ and (iv) $L$ intersects
$\p (L)$ transversally. Then there is a natural isomorphism between 
\[I^{(r_0)}_*(L, \p_s; P) \cong H^{* + m}(L; Z_2)  \ \ \ \ 
\mbox{for $* \in Z$ and $s \in [0,1]$}. \]
\end{pro}
Proof: For any $s, s^{'} \in [0,1]$ sufficiently close, 
by Theorem~\ref{invariant} we have  
\begin{equation} \label{ss}
I^{(r_0)}_*(L, \p_s; P) \cong I^{(r_0)}_*(L, {\p }_{s^{'}}; P) .
\end{equation}
For $s \in [0,1]$ sufficiently small, we have
\begin{equation} \label{so}
I^{(r_0)}_*(L, \p_s; P) \cong H^*(C_s, \bd_s) \cong H^{*+m}(L; Z_2).
\end{equation}
The first isomorphism is given by Lemma~\ref{cs} and the second by Lemma 3 in 
\cite{ch}. Then the result follows from finite steps of applying
(\ref{ss}) and (\ref{so}). \qed
\medskip

\noindent{\bf Remarks}: (i) 
Proposition~\ref{chek} provides the Arnold conjecture
for monotone Lagrangian sub-manifold with $\Si (L) \geq 3$ 
and $E(\p ) < \s (L)$. Chekanov's result in \cite{ch} 
does not require the assumptions of the 
monotonicity of $L$ and $\Si (L) \geq 3$. Floer proved the Arnold conjecture
for monotone case. The generalized version has been obtained by Hofer and Salamon
\cite{hs}, and Ono \cite{on} for weakly monotone case.
As in \cite{oh1} Remark 3.7, 
we can also use the integer graded Floer cohomology to
show Gromov's non-exactness theorem for compact Lagrangian embeddings into
$C^n$.

(ii) Using
$I^{(r_0)}_*(L, \p_s; P)$, there is a natural relation 
$I^{(r_0 + \s (L))}_*(L, \p_s; P) = I^{(r_0)}_{* + \Si (L)}(L, \p_s; P)$ for
different choices of $a_s, a_s + \s (L)$. In fact 
\[I^{(r_0)}_*(L, \p_s; P) \cong I^{(r_0)}_*(L, {\p }_{s^{'}}; P), \]
for $s, s^{'} \in [0,1]$ which answers Oh's question in \cite{oh2} page 29.
\medskip

For Lagrangian sub-manifold $L \subset (C^m, \o_0)$ with the standard symplectic
structure $\o_0 = - d \lam $, where $\lam $ is the Liouville form, there is a
Liouville class $[\lam |_L] \in H^1(L, \R)$. One of the fundamental results
in \cite{gr} is the non-triviality of the Liouville class.
 
\begin{thm}[Gromov \cite{gr}] \label{gromov} For any compact Lagrangian
embedding $L$ in $C^m$, the Liouville class $[\lam |_L] \neq 0 \in H^1(L, \R)$.
\end{thm}

In particular, $H^1(L, Z_2) \neq 0$.

\begin{thm} \label{raudin}
For any small $\ve > 0$, if $e_H(L) = \s (L) - \ve$ for a 
monotone Lagrangian sub-manfold $L$ embedded in $C^m$,
then $\Si (L) \leq 2$.
\end{thm}
Proof: Suppose the contrary, $\Si (L) \geq 3$. By the hypothesis, there exists a
$\p \in {\cal D}_{\o}$ such that 
\begin{equation}
 e_H(L) < E(\p ) \leq e_H(L) + \frac{\ve }{2} < \s (L).
\end{equation} 
So the time one flow $\p$ separates the $L$ from its definition of $e_H(L)$, 
i.e. $L \cap \p (L) = \emptyset $. By definition of integer graded
symplectic Floer cohomology, we have $I^{(r_0)}_*(L, \p; P) = 0$.
By Proposition~\ref{chek}, we have
\[ I^{(r_0)}_*(L, \p_s; P)\cong H^{* + m}(L; Z_2) \ \ \ \ \  
\mbox{for $s \in [0,1]$}.\]
This contradicts with Theorem~\ref{gromov} 
$H^*(L; Z_2) \neq 0$. So we obtain the restriction on the Maslov index. \qed

\subsection{Lagrangian rigidity and Audin's question}

\begin{pro} For any compact monotone Lagrangian sub-manifold in $(P, \o)$, if
$\Si (L) \geq m +1, (m \geq 2)$, then
\begin{enumerate}
\item all the differentials $d^k$ are trivial for $k \geq 0$,
\item we have the following relations.
\[ \sum_{k \in Z}I^{(r_0)}_{j + \Si (L) k } (L, \p; P) = HF^j(L, \p; P) . \]
\end{enumerate}
\end{pro}
Proof: By Proposition~\ref{chek}, for $s$ sufficiently 
small the Maslov index for $I^{(r_0)}_*(L, \p_s;P,J)$ satisfies the following.
\[0 < \max \mu^{(r_0)}(y^{(r_0)}) - \min \mu^{(r_0)}(x^{(r_0)}) =
\mu_{H_s}(y) - \mu_{H_s}(x) \leq m. \]
The result follows from the definition of $d^k$. \qed

\begin{thm}[Oh \cite{oh1}] \label{rigi}
For any compact monotone Lagrangian embedding $L \subset C^m$,
\[ 1 \leq \Si (L) \leq m. \]
\end{thm}
 
Note that the result is optimal based on examples in \cite{po2}. We can use
the integer graded symplectic Floer cohomology to prove Theorem~\ref{rigi} by
the same argument used for the Floer-Oh's local Floer cohomology. 
We define the associated Poincar\'{e}-Laurent polynomials (shifted in degree) 
for the
spectral sequence $P(E^k, t), k \geq 1$ by
\begin{equation} \label{plp}
P^{(r)}(E^k, t) = \sum_{n \in Z} (\dim_{Z_2} E^k_{n,j}) t^{n}. 
\end{equation}

 Note that our polynomials are slightly different from the one formulated in
\cite{oh1}.
By Theorem~\ref{E1}, $P^{(r)}(E^1, t) = \sum_{n \in Z} 
(\dim_{Z_2} I^{(r)}_n(L, \p ; P, J)) t^{n}$. 
From \S 5.1 Remarks (ii), we have
\begin{equation} \label{shift}
P^{(r+\s (L))}(E^k, t) t^{\Si (L)} = P^{(r)}(E^k, t).
\end{equation}

\begin{pro} \label{pl} For monotone Lagrangian $L$ with $\Si (L) \geq 3$,
\[P^{(r)}(E^k, t) = \sum^k_{i =1} (1 + t^{-i \Si (L) -1}) \overline{Q}_i(t) +
P^{(r)}(HF^*, t), \]
where $k+1 = k(L)$ and $\overline{Q}_i(t)$ 
are Poincar\'{e}-Laurent polynomials of
nonnegative integer coefficients.
\end{pro}
Proof: Let $Z^1_{n,j} = \ker \{d^1: E^1_{n,j} \to E^1_{n+ \Si (L) +1, j+1}\}$
and $B^1_{n,j} = Im d^1 \cap E^1_{n,j}$. we have the exact sequences
\[0 \to Z^1_{n,j} \to E^1_{n,j} \to B^1_{n+ \Si (L) +1, j+1} \to 0, \]
\[0 \to B^1_{n,j} \to Z^1_{n,j} \to E^2_{n,j}  \to 0. \]
So the degree $\Si (L) +1$ of $d^1$ derives the followings.
\begin{equation} \label{e12}
P^{(r)}(E^1, t) = P^{(r)}(E^2, t) + (1+t^{- \Si (L) - 1}) P^{(r)}(B^1, t).
\end{equation}
Since the higher differential $d^i$ has degree $i \Si (L) +1$, we can repeat
(\ref{e12}) and let $\overline{Q}_i(t) = P^{(r)}(B^i, t)$.
Note that $E^{\infty } \cong HF^*(L, \p; P)$ by Theorem~\ref{invariant}. Thus
we obtain the desired result. \qed

For any
oriented monotone Lagrangian torus in $C^m$, its Maslov number is always
even. Suppose $\Si (L) \neq 2$. Then $\Si (L) \geq 4 ( > 3)$. 
From the construction of integer graded symplectic Floer cohomology,
$I^{(r_0)}_*(L, \p; P)$ is well-defined 
provided $\Si (L) \geq 3$. Also the Floer cohomology $HF^*(L, \p; P)$
is well-defined and is invariant under the generic continuation of
$(J, \p )$ from \cite{oh}. 
By choosing a Hamiltonian isotope $\p = \{ \p_t\}$
such that
\[ L \cap \p_1 (L) = \emptyset , \]
which is certainly possible in $C^m$. Thus $Z_{{\p }_1} = \emptyset$, and
\begin{equation} \label{trivial}
HF^*(L, \p; P) = 0 .
\end{equation}
On the other hand, while $\p_s $ is sufficiently close to identity map with 
$E(\p_s ) < \s (L)$, by Proposition~\ref{chek} we have
\[I^{(r_0)}_*(L, {\p }_s; P) \cong H^{*+m}(L; Z_2). \]
Applying Proposition~\ref{pl} for $E^1 = H^{*+m}(L; Z_2)$, 
$E^{\infty} = HF^*(L, \p;
P) = 0$ from (\ref{trivial}), we get the relation
\begin{equation} \label{pt}
(1+t)^m t^{-m}= \sum_{i=1}^k (1 + t^{- i \Si (L) - 1}) \overline{Q}_i(t),
\end{equation} 
where $(1+t)^mt^{-m}$ is the Poincar\'{e}-Laurent polynomial for 
the torus $T^m$.
Note that our Poincar\'{e}-Laurent polynomial (\ref{pt}) is similar to
the one
in \cite{oh1}, \cite{oh2}. The 
changing $+1$ into $-1$ reflects precisely 
the degrees $\Si (L) + 1, - \Si (L) + 1$
of differentials in the spectral sequences 
for the Floer cohomology. After shifting the degree, we just have
the following Poincar\'{e} polynomial:
\begin{equation} \label{ohp}
(1+t)^m = \sum_{i=1}^k (1 + t^{i \Si (L)+1}) Q_i(t),
\end{equation}
where $Q_i(t) = t^{m- i \Si (L) - 1}\overline{Q}_i(t)$ is a polynomial
in $t$ with nonnegative integer coefficients.

\begin{lm} For an embedded, oriented, monotone Lagrangian
torus $L$ in $C^m$ with $k(L) \leq 2$, we have
\[ \Si (L) = 2. \]
\end{lm}
Proof:
(i) If $k(L) = 1$, then $k =0$, we get a contradiction of $HF^*(L, \p; P) = 0$
and $E^{\infty } = I^{(r_0)}_*(L, \p; P) \cong H^*(L; Z_2) \neq 0$ by 
Theorem~\ref{gromov}.

(ii) If $k(L) =2$, then by Proposition~\ref{pl} and (\ref{ohp}) 
we have the identity
\[ (1+t)^m = (1+t^{ \Si (L) +1})Q_1(t) .\]
But this is impossible to have
such a decomposition of $(1+t)^m = (1+t^{\Si (L)+1})q_1(t)$ for any even 
$\Si (L) \geq 3$.
So $\Si (L)$ has to be an even integer in
$2 \leq \Si (L) < 3$, i.e. $\Si (L) = 2$. \qed

\noindent{\bf Remark:}
Note that the result can be restated as that for any given compatible almost
complex structure $J$, the monotone Lagrangian torus $L$ carry a $J$-holomorphic
disc $u: (D^2, \bd D^2) \to (C^m, L)$ with $\mu_L(u) = 2$. Following the
discussion in \cite{oh1}, one has $0 < k \Si (L) \leq m+1$, i.e
\[1 \leq k(L) = k + 1 \leq [\frac{m+1}{\Si (L)}] +1. \]
In particular, if $\Si (L) \geq 3$ and $\Si (L) | (m+1)$, then the last term
in (\ref{ohp}) may be
\[ (1 + t^{[\frac{m+1}{\Si (L)}] \Si (L) + 1} Q_{[\frac{m+1}{\Si (L)}]}(t) 
= (1 + t^{m+2}) Q_{[\frac{m+1}{\Si (L)}]}(t). \]
Comparing with the left hand side of (\ref{ohp}), we have $
Q_{[\frac{m+1}{\Si (L)}]}(t) = 0$, i.e. $k(L) = k +1 \leq (m+1)/{\Si (L)}$
in this case.

\begin{thm} \label{unless}
For an embedded, oriented, monotone Lagrangian
torus $L$ in $C^m$, we have
\[ \Si (L) = 2, \]
unless $k(L) = \frac{m+1}{\Si (L)}$.
\end{thm} 
Proof: Suppose the contrary. So $\Si (L) \geq 3$. Our integer
graded symplectic Floer cohomology is well-defined for the
transverse intersection $L \cap \p (L)$. We have
$I^{(r_0)}_*(L, \p_s;P) \cong H^{*+m}(L;Z_2)$. So 
$E^1_{n,j}$ is contributed from the cohomology of
torus with corresponding degree shift, i.e.
\[H^j(T^m;Z_2) = \oplus_{n \equiv j mod \Si (L)} E^1_{n,j} .\]
Note that our higher differentials
\[ d^k: E^k_{n,j} \to E^k_{n+k \Si (L) +1, j+1},\]
from $H^j(T^m;Z_2)$ to $H^{j+1}(T^m;Z_2)$ modulo $E^{k-1}$. By (\ref{trivial}),
$HF^*(T^m, \p; P, Z_2) = 0$, so there are no elements in $E^*$ survived from
the differentials. By counting the rank, we have
\begin{equation} \label{sum0}
\sum_{j=0}^{\Si (L) -1} \sum_{i=1}^k (-1)^j 
rank(\oplus_{n \equiv j mod \Si (L)} E^i_{n,j}) = 0,
\end{equation}
where $j = 0, 1, \cdots, \Si (L) -1$.
Clearly we have, for each $j = 0, 1, \cdots, \Si (L) -1$,  
\[\sum_{i=1}^k rank(\oplus_{n \equiv j mod \Si (L)} E^i_{n,j}) =  
rank H^{j}(T^m ; Z_2) \]\[
+ rank H^{j+ \Si (L)}(T^m ; Z_2) + \cdots
+ rank H^{j+ k \Si (L)}(T^m ; Z_2).\]
Thus (\ref{sum0}) gives us the following identity
\begin{equation} \label{sum1}
\sum_{l=0}^k  \left( \begin{array}{c} m \\ l \Si (L) \end{array}
\right) - 
\sum_{l=0}^k  \left( \begin{array}{c} m \\ l \Si (L) + 1 \end{array}
\right) + \cdots - \sum_{l=0}^k  \left( 
\begin{array}{c} m \\ l \Si (L) + \Si (L) -1 \end{array}
\right) = 0 .
\end{equation}
This can be rearrange into the alternating sum in the following 
due to the even number $\Si (L)$.
\begin{equation} \label{sum2}
\sum_{l=0}^{(k+1) \Si (L) - 1}(-1)^l \left( \begin{array}{c} m \\ l \end{array}
\right) = 0. 
\end{equation}
This derives another constraint 
\begin{equation} \label{sum3}
(k+1) \Si (L) - 1 \leq m,
\end{equation}
i.e. $(k+1) \leq \frac{m+1}{\Si (L)}$ which is better restriction than 
$(k+1) \leq \frac{m+1}{\Si (L)}+1$. So we have $(k+1) \leq [\frac{m+1}{\Si (L)}]$.

(i) If $(k+1) \leq [\frac{m+1}{\Si (L)}] - 1$, then
\begin{eqnarray*}
(k+1) \Si (L) - 1 & \leq & ([\frac{m+1}{\Si (L)}]- 1)\Si (L) - 1 \\
 & \leq & m - \Si (L),
\end{eqnarray*}
for $\Si (L) \geq 3$.
Thus the equation (\ref{sum2}) will be equivalent to
\begin{equation} \label{sum4}
\sum_{l= (k+1) \Si (L)}^m(-1)^l \left( \begin{array}{c} m \\ l \end{array}
\right) = 0 ,
\end{equation}
by $(1 -1)^m = 0$. This is impossible for $m \geq 3$ and (\ref{sum4}) with
$4 \leq (k+1) \Si (L) < m-2$. So the contradiction shows that
$\Si (L) = 2$ for this case.

(ii) If $k = [\frac{m+1}{\Si (L)}]$, $m+1 = k \Si (L) + r$,
$0 \leq r \leq \Si (L) -1$, then by the discussion above, we have
$k \leq \frac{m+1}{\Si (L)} -1$ if $r =0$. So we only consider
$1 \leq r \leq \Si (L) -1$.
In this case, we obtain
\[(k+1) \Si (L) - 1 = m -r + \Si (L) \geq m + 1,\]
which contradicts with (\ref{sum3}).

(iii) If $(k+1) = [\frac{m+1}{\Si (L)}]$ for $r \neq 0$ in
\[ m+1 = (k+1) \Si (L) +r . \]
It is easy to see that 
\[(k+1) \Si (L) - 1 = m -r < m .\]
So by the same argument in (i), we have a contradiction. 
Therefore the result follows except the case $m+1 = (k+1) \Si (L)$ (i.e.
$k(L) = \frac{m+1}{\Si (L)}$). \qed

Since $\Si (L)$ is always even for an oriented Lagrangian sub-manifold $L$,
we offered an affirmative answer to the Audin's question for all even dimension
$m$ by Theorem~\ref{unless}. Now we are going to complete the proof of Theorem
C by considering the odd dimension $m$ ($m \equiv 1 \pmod 2$) with the
divisibility $\Si (L) | (m +1)$.

For an embedded, oriented, monotone Lagrangian torus $L$ in $(C^m, \o)$, 
we have $I_{\o } = \lam I_{\mu, L}$ for $\lam \geq 0$. It is not
hard to see that $L \x L$ is monotone in $(C^m \x C^m, \o \oplus \o)$.
For $u: (D^2, \bd D^2) \to (C^m \x C^m, L \x L)$, we denote
$p_i: C^m \x C^m \to C^m$ to be the projection on the $i$-th factor,
$i=1, 2$. Thus one has
\[I_{\mu \oplus \mu , L \x L}(u) = I_{\mu, L}(p_1 u) + I_{\mu, L}(p_2 u), \]
which follows from the product symplectic form and the
K\"{u}nneth formula for the Maslov class in $\Lam (C^m) \x \Lam (C^m)$
(see \cite {clm}). Hence we have
\begin{eqnarray*}
I_{\o \oplus \o}(u) & = & I_{\o } (p_1 u ) + I_{\o }(p_2 (u) \\
 & = & \lam I_{\mu, L}(p_1 u) + \lam I_{\mu, L}(p_2 u) \\
 & = & \lam I_{\mu \oplus \mu , L \x L}(u) .
\end{eqnarray*}

\begin{pro} \label{odd}
For $m$ odd, $\Si (L) | (m+1)$, $L$ is an oriented monotone Lagrangian
embedding into $C^m$, we have $\Si (L) = 2$.
\end{pro}
Proof: Suppose the contrary $\Si (L) > 3$. Combining the discussion above, we
have a monotone, oriented embedding
\[ L \x L \hookrightarrow (C^m \x C^m, \o \oplus \o) .\]
Note that $\Si (L \x L) \geq \Si (L)$ from the additivity of the Maslov index,
we have $\Si (L \x L) > 3$, and $L \x L$ is an oriented
monotone Lagrangian embedding into $2m$ {\bf even} dimensional
complex plane. Thus by Theorem~\ref{unless}, we have
$2 = \Si (L \x L) \geq \Si (L) > 3$ which is absurd.
Thus $\Si (L) =2$ for the last case. \qed

\noindent{\bf Remark:} Using the same method, we can also obtain the result
from the Floer-Oh's local homology and the 
spectral sequence defined in \cite{oh1}. 
The point is to use the full spectral sequence information, rather
than just the identity from the Pincar\'{e} polynomial of the spectral sequence.
\medskip

There are four interesting numbers of a monotone Lagrangian manifold $L$, 
$\Si (L), \s (L), e_H(L), k(L)$. Both of them intertwines together and links
with the integer graded symplectic Floer cohomology and its derived spectral
sequence. It is interesting to study further relations among them.
\medskip

\noindent{\bf Acknowledgment}: I would like to thank 
Y.G. Oh for his explanation on his work and many helps.
It is also pleasure to thank I. Newton Institute and Max-Planck Institute for 
their financial support through the visit and for nice research environment.
We would like to thank the organizers S. Donaldson, D. McDuff, D.Salamon and
C. Thomas for making it possible to participate a part of the program 
``Symplectic Geometry'' in the fall of 1994.

\end{document}